\documentclass[11pt]{elsarticle}
\usepackage[utf8]{inputenc}
\usepackage[T1]{fontenc}
\usepackage{amsfonts}
\usepackage{mathrsfs,amsmath}
\usepackage{algorithm, algorithmic}
\usepackage{color}
\usepackage{bm}
\usepackage{hyperref}
\usepackage{fullpage}
\usepackage{empheq}
\usepackage{multirow}
\usepackage{subcaption}
\usepackage{booktabs,array}
\usepackage{multicol}
\usepackage{float}
\usepackage{amsthm}
\usepackage{here}
\usepackage{siunitx}

\theoremstyle{definition}
\newtheorem{remark}{Remark}

\DeclareMathOperator{\sign}{sign}

\begin{document}

\title{Kinetic scrape off layer simulations with semi-Lagrangian discontinuous Galerkin schemes}
\author[uibk]{Lukas Einkemmer\corref{cor1}} \ead{lukas.einkemmer@uibk.ac.at}
\author[uibk]{Alexander Moriggl}
\address[uibk]{Department of Mathematics, Universit\"at Innsbruck, Austria}
\cortext[cor1]{Corresponding author}

\begin{abstract} In this paper we propose a semi-Lagrangian discontinuous Galerkin solver for the simulation of the scrape off layer for an electron-ion plasma. We use a time adaptive velocity space to deal with fast particles leaving the computational domain, a block structured mesh to resolve the sharp gradient in the plasma sheath, and limiters to avoid oscillations in the density function. In particular, we propose a limiter that can be computed directly from the information used in the semi-Lagrangian discontinuous Galerkin advection step. This limiter is particularly efficient on graphic processing units (GPUs) and compares favorable with limiters from the literature. We provide numerical results for a set of benchmark problems and compare different limiting strategies.
\end{abstract}


\maketitle
 
\section{Introduction}

The scrape of layer (SOL) in a fusion reactor is characterized by open field lines that connect the plasma to the divertor plates. In this region the plasma is no longer confined and thus high energy particles hit the solid material at the divertor. The corresponding heat flux is an important consideration for reactor design (see, e.g, \cite{Freidberg2015,Myra2011,Coulette2016,Gasteiger2017,Bourne2023}). A number of interesting physical phenomena occur in the SOL. For example, fast particles travel along magnetic field lines to the wall and are lost. This happens more quickly for electrons (which have much larger average speed) than for ions, resulting in an electric field that accelerates the ions. Close to the wall, a plasma sheath forms which is characterized by a sharp gradient in the electric potential. Such phenomena are difficult to capture with numerical methods.

The fundamental model to describe such a problem are the Vlasov--Poisson equations for the electron and ion density. The two particle species are coupled through the electric field. There are a number of methods to solve this system numerically. For example, particle in cell schemes (see, e.g. \cite{Verboncoeur2005}), Eulerian methods (see, e.g., \cite{Arber2002, Filbet2003, Heath2012, Hittinger2013}), semi-Lagrangian methods (see, e.g. \cite{Cheng1976,Klimas1994,Sonnendruecker1999,Crouseilles2011,Rossmanith2011,Qiu2011,Crouseilles2009,Einkemmer2019b,AllmannRahn2022, Guo2022,Guo2022b, Bourne2023}), and various types of complexity reduction schemes (see, e.g, \cite{Kormann2014,Kormann2015,Einkemmer2018,Einkemmer2020,Einkemmer2021a, Cassini2021,Coughlin2022,Coughlin2023,Einkemmer2023}) have been proposed. In this paper we focus on semi-Lagrangian schemes. These schemes follow the characteristics backward in time and perform an interpolation/projection to obtain the density at the new time. Such methods are attractive as they do not suffer from a Courant--Friedrichs--Lewy (CFL) condition. Thus, the time step size is independent of the largest particle speed and the cell size. There is a variety of semi-Lagrangian approaches in the literature. Lagrange interpolation is simple, but introduces excessive numerical diffusion. Using splines performs better in this regard and cubic splines are commonly used (see, e.g. \cite{Filbet2003,Grandgirard2006a,Bourne2023}). Spectral methods can also be considered semi-Lagrangian schemes. Since such methods tend to oscillate and thus filtering is required \cite{Klimas1994}. Both spline interpolation and spectral methods have global data dependency and thus are problematic for large scale parallel computing. Although more local variants of cubic spline interpolation have been proposed \cite{Crouseilles2009}, those sacrifice some of the properties of the scheme. We also note that most of the recent literature on semi-Lagrangian schemes on graphic processing units (GPUs) uses Lagrange interpolation and not splines (see, e.g., \cite{Mehrenberger2013,Kormann2019}).

More recently, semi-Lagrangian discontinuous Galerkin schemes (sLdG) have been proposed \cite{Crouseilles2011,Rossmanith2011,Qiu2011}. These have the advantage that they have local data dependency and comparable or in some cases even superior \cite{Einkemmer2019b} properties compared to cubic spline based semi-Lagrangian methods. Thus, they are ideally suited for large scale parallel computations including GPUs \cite{Einkemmer2020a,Einkemmer2024}. In \cite{Einkemmer2021c} simulations of six-dimensional problems have been conducted on supercomputers with thousands of GPUs. In this paper we address the challenges that are posed by SOL simulation for such methods. In fact, the local nature of the sLdG method turns out to be very beneficial in this case. In particular, we resolve the sharp gradient in the plasma sheath by a fine mesh that couples via an (artificial) boundary to the coarser mesh used in the remainder of the domain. Since the discontinuous Galerkin approximation uses polynomials in each cell transferring them from a coarse mesh to a fine mesh is immediate and the reverse can be done by a simple projection. We also use a time adaptive cut-off for the velocity space to efficiently deal with the fact that fast particles leave the computational domain early in the numerical simulation. 

Finally, we consider limiters to suppress oscillations in the simulation. While it is true that in many cases Vlasov--Poisson simulations can be conducted without limiters. In fact, any limiting introduces additional numerical diffusion which, in general, is unwanted and can interfere with the preservation of some physical invariants \cite{Crouseilles2010,Einkemmer2017a}. However, for the SOL the sharp gradient in the plasma sheath make limiters a necessity. Many such approaches have been proposed in the literature (see, e.g. \cite{Crouseilles2010,Rossmanith2011,Qiu2011,Guo2013}). An advantage of the sLdG approach is that limiters are generally applied locally and ideas and methods from the vast literature on Runge--Kutta discontinuous Galerkin methods (see, e.g., the reviews \cite{Cockburn2000,Qiu2005,Kochi2023}) can be used. However, one downside of such limiters is that they are formulated as a post processing step requiring data from the cells to the left and to the right of the current cell. This can add significant computational cost for a memory bound problem, where the direction of flow is known. Thus, we propose a new limiting scheme that can be directly applied in the sLdG time step and thus has significantly lower computational cost. 

The remainder of the paper is structured as follows. In section \ref{sec:model} we introduce the Vlasov--Poisson equations for electrons and ions. Then in section \ref{sec:sldg} the time splitting based sLdG scheme is introduced. In section \ref{sec:limiter} we describe the limiters and discuss their characteristics. In section \ref{sec:amr} we present the adaptive mesh refinement scheme (which consists of an adaptive adjustment of the velocity domain and a block structured refinement close to the wall). Finally, we present numerical results in section \ref{sec:numerical}.

\section{Model equations \label{sec:model}}

We consider an electrostatic kinetic model for each species of the plasma. More specifically, we have (see, e.g., \cite{Coulette2016})
\begin{equation}
\begin{aligned}
 \partial_t f_\alpha(t,x,v) + v\cdot \partial_x f_\alpha(t,x,v) + \frac{q_\alpha E(t,x)}{m_\alpha} \partial_v f_\alpha(t,x,v) &= S_\alpha(t,x,v)  \\
 -\partial_{xx} \phi(t,x) = \sum_\alpha q_\alpha \int f_\alpha(t,x,v) \,dv, \quad E(t,x) &= -\partial_x \phi(t,x)
\end{aligned}
\end{equation}
where $\alpha$ is the particle species (e.g.~electrons or ions). The charge and mass of particle species $\alpha$ are denoted by $q_\alpha$ and $m_\alpha$, respectively, and $S_\alpha$ is a source term. We will restrict ourselves to the two-dimensional (one dimension in space and one dimension in velocity) case here. The spatial domain is $x \in [-L, L]$.  In addition, to the kinetic equation a Poisson problem has to be solved to obtain the electric potential. The electric field $E$ is then determined from this potential. 

 We can also include a collision operator on the right-hand side. In fact, this has been implemented in our code. However, in this work we will focus on the collisionless case. This is the most difficult case with respect to limiting and we emphasize that realistic collision rates are not sufficient to suppress oscillations. In fact, numerical simulations have shown that compared to the physical value the collision rate has to be increased by a factor of $100$-$1000$ in order to obtain reasonable results.

For SOL simulations it is important to impose appropriate boundary conditions. Let us consider this at the right boundary, i.e.~at $x=L$. Then for $v>0$ we have particles that hit the wall. We impose an outflow condition in this case (i.e.~particles leave the computational domain). In addition, we impose zero inflow. That is, we set $f(t,L,v)=0$ for $v<0$. For the potential $\phi$ homogeneous Dirichlet boundary conditions are assumed. That is, a constant potential at the divertor plates is imposed.

In this paper, we consider two particle species, electrons and ions. After nondimensionalizing the system according to table~\ref{tab:nondimensionalizing} (see, e.g., \cite{Bourne2023} for more details), we end up with the following equations
\begin{equation}
\begin{aligned}
 \partial_t f_e + v \partial_x f_e - E(x) \partial_v f_e = S, \\
 \partial_t f_i + \frac{v}{\sqrt{\mu}} \partial_x f_i + \frac{E(x)}{\sqrt{\mu}} \partial_v f_i = S, \\
 -\partial_{xx} \phi = \int f_i\,dv - \int f_e\,dv, \quad E(x) = -\partial_x \phi,
\end{aligned}
\end{equation}
where $\mu = \frac{m_i}{m_e}$ is the mass ratio between ions and electrons. Note that the velocity normalization for electrons and ions are different. For both particle species the same source term is used (i.e.~the source is neutral).

\begin{table}
\centering
\begin{tabular}{l|lcl}
 space:    & $\hat{x}$ & = & $ x/\lambda$\\ 
 velocity: & $\hat{v}$ & = & $v/v_{\text{th},\alpha}$ \\ 
 time:     & $\hat{t}$ & = & $t  v_{\text{th},e}/\lambda$\\
 particle density:  & $\hat{\rho}_\alpha$ & = & $ \rho_\alpha/\rho_0 $ \\
 distribution density:  & $\hat{f}_\alpha$ & = & $f_\alpha  v_{\text{th},\alpha}/\rho_0$\\
 electric field: & $\hat{E}$ & = & $E q_e \lambda / T$ \\
 electric potential & $\hat{\phi}$ & = & $\phi  q_e/T$ \\
 \end{tabular}
 \caption{Nondimensionalization of the Vlasov--Poisson equations. Spatial variables are normalized with respect to the Debye length $\lambda$, the velocity is normalized with respect to the thermal velocity of electrons and ions, respectively, number density is normalized with respect to the initial background number density $\rho_0$. The electric field is normalized such that there there are no dimensionless constants in Gauss's law. All other normalizations follow from this choice. Note that we use $T = m_e v_e^2$ in the table.  \label{tab:nondimensionalizing}}
\end{table}

\section{Time splitting/semi-Lagrangian discontinuous Galerkin method} \label{sec:sldg}

We start by performing a splitting scheme for the Vlasov equation. Such an approach has been commonly used in the literature \cite{Cheng1976,Sonnendruecker1999,Filbet2003}. It has the advantage that we only need to solve one-dimensional advections, which will be done using the sLdG scheme, and a set of ordinary differential equations for the source term. Moreover, the properties of these methods and their convergence has been studied extensively in the mathematics literature \cite{Besse2008,Einkemmer2014b,Einkemmer2014a}. The second order accurate Strang splitting scheme is shown in Algorithm \ref{alg:splitting_1}.

\begin{algorithm}
\caption{Second order Strang splitting scheme. We use $\mu_{\alpha} = m_{\alpha}/m_e$ and the charge of the particles is denoted by $c_\alpha$, i.e.~$c_e=-1$ and $c_i=1$).}
\label{alg:splitting_1}
\begin{enumerate}
\item \textbf{Source:} Perform a RK2 step for $\partial_t f_\alpha = S_\alpha $ with step size~$\Delta t/2$ and initial condition $f^n_\alpha$ to get $f^*_\alpha$.  \\
\item \textbf{Advection x:} Perform a sLdG step for the one-dimensional advection $\partial_t f_\alpha + \frac{v}{\sqrt{\mu_\alpha}}\partial_xf_\alpha=0$ with step size $\Delta t/2$ and initial condition $f_\alpha^*$ to get $f_\alpha^{**}$.
\item \textbf{Poisson:} Solve the Poisson problem $-\partial_{xx} \phi = \sum_\alpha q_\alpha \int f_\alpha^{\star \star}\,dv$  and use it to compute $E^{n+1/2}$. 
\item \textbf{Advection v:} Perform a sLdG step for the one-dimensional advection $\partial_tf_\alpha + \frac{c_\alpha E^{n+1/2}}{\sqrt{\mu_\alpha}}\partial_v f_\alpha =0$ with step size~$\Delta t$ and initial condition $f_\alpha^{**}$ in order to obtain $f_\alpha^{***}$,
\item \textbf{Advection x:} Perform a sLdG step for the one-dimensional advection $\partial_t f + \frac{v}{\sqrt{\mu_\alpha}}\partial_x f_\alpha=0$ with step size~$\Delta t/2$ and initial condition $f_\alpha^{***}$ to obtain $f_\alpha^{****}$, 
\item \textbf{Source:} Perform a RK2 step for $\partial_t f_\alpha = S_\alpha$ with step size~$\Delta t/2$ and initial condition $f^{****}_\alpha$ to obtain the solution at the new time-step $f^{n+1}_\alpha$,
\end{enumerate}
\end{algorithm}

Note that since step 4 does not change $\int f_{\alpha} \,dv$ the value computed in step 3 of Algorithm \ref{alg:splitting_1} indeed is an approximation to $E(t_n+\Delta t/2)$. The approximation is first order accurate which is sufficient to get a second order accurate scheme overall (see \cite{Sonnendruecker1999} for more details).

\subsection{Advection step}

The splitting schemes reduces the Vlasov--Poisson equation to a set of one-dimensional advection equations. In order to solve these equations, we employ a semi-Lagrangian discontinuous Galerkin scheme. Both the advection in $x$ and $v$ can be written in the following form
\begin{equation} \label{eq:advection1d}
\left\{
\begin{aligned}
 \partial_t u(t,x) + a\partial_x u(t,x) &= 0 \\
 u(0,x) &= u_0(x),
 \end{aligned}
 \right.
\end{equation}
where $a$ is the advection speed. Since the advection speed for the $x$ advection depends only on $v$ and vice versa, $a$ can be considered constant for the purpose of solving the advection equation. 

Equation \eqref{eq:advection1d} has the solution $u(t,x) = u_0(x-at)$, i.e., the solution is nothing else than a shift of the initial condition. In the discontinuous Galerkin framework, $u_0(x)$ is a piecewise function, which is discontinuous at the cell interfaces and polynomial of degree $k$ inside each cell. In our implementation, Lagrange polynomials which interpolate at the Gauss--Legendre nodes are used to represent the polynomials. Thus, $k+1$ function values at the Gauss--Legendre nodes are stored in computer memory for each cell. Shifting this piecewise polynomial representation by $-a \Delta$, in general, shifts the discontinuity to a location inside the cells. In order to obtain a polynomial approximation in each cell, the discontinuous function is projected back to the space of polynomials. This $L^2$ projection is of order $k+1$ and conservative, see~\cite{Crouseilles2011,Rossmanith2011,Qiu2011}. Its coefficient can be computed analytically which yields
\begin{equation}
\label{eq:sldg_implementation}
 u^{n+1}_{i,\cdot} = A(\alpha)u^n_{i^*,\cdot} + B(\alpha)u^n_{i^*+1,\cdot}
\end{equation}
where $-a\Delta t = \Delta x \cdot (i^* + \alpha)$, $\alpha \in [0,1)$ and $i^* \in \mathbb{Z}$. Thus, the solution at the new time-step is simply the sum of two small matrix-vector products of two adjacent cells. The matrices $A(\alpha)$ and $B(\alpha)$ are precomputed and can then used for every cell in the same one-dimensional slice of the domain.

This sLdG approach has been considered extensively in the literature. We have followed here the notation from \cite{Einkemmer2016} and we thus refer the reader there for more details. Let us remark that this method does not suffer from any stability restriction and thus large step sizes can be taken independent of the order of the method and the grid spacing.

\subsection{Poisson problem}

In order to solve the one-dimensional Poisson equation, we use a symmetric interior penalization Galerkin (SIPG) scheme (see, e.g., \cite{Riviere2008}). This scheme is local, symmetric, and homogeneous Dirichlet boundary conditions can be easily imposed. Moreover, a domain with different cell sizes can be used in a straightforward manner. This also has the advantage that we can directly use the high-order discontinuous Galerkin approximation space as for the advection problem.

However, in order to obtain the electric field, the solution of the Poisson equation, i.e., the electric potential, has to be differentiated. This reduces the order of the method by one, which was also observed in~\cite{Cai2019}. To overcome this, we first evaluate the right-hand (polynomials of degree $k$) to the grid which is defined by polynomials of degree $k+1$. Then, the electric potential is obtained with order $k+2$ which gives an electric field of order $k+1$.

Since the Poisson problem is only one-dimensional, it is comparatively cheap to solve. The linear system that arises from the SIPG method is symmetric and block-tridiagonal. Therefore we can factorize the matrix once and for all at the beginning of the simulation. The factorization can than be reused at every time-step to solve the Poisson equation quickly. We have used the general banded solver from the fast Lapack library~\cite{lapack} to perform this task. Let us note that for higher-dimensional problem a similar approach can be used, but it is more efficient to solve the linear system using e.g.~a conjugate gradient method (see, e.g., \cite{Einkemmer2021c}).

\section{Limiters \label{sec:limiter}}

In many application of semi-Lagrangian discontinuous Galerkin schemes for kinetic equations no limiting is necessary, i.e.~the small but inherent numerical diffusion in the scheme is sufficient to obtain sufficiently accurate solutions \cite{Einkemmer2019b}. In fact, as has been observed in \cite{Crouseilles2010,Einkemmer2017a} one has to be careful as adding limiters (e.g.~to enforce positivity) can add additional numerical diffusion which, in principle, is unwanted and can degrade the numerical solution. However, due to the sharp gradients in the plasma sheath some amount of limiting is necessary to suppress oscillations, in particular, in the electron density function. 

In the early works of semi-Lagrangian discontinuous Galerkin schemes limiters were introduced with the goal of obtaining a positive solution \cite{Rossmanith2011,Qiu2011}. In addition, limiters for Runge--Kutta discontinuous Galerkin (RKDG) schemes, which are widely used for simulating conservation laws, have been developed. These schemes are related to sLdG schemes as they share the same spatial approximation; see, e.g., \cite{Cockburn2000,Cockburn2001}. Such limiters work as follows. First, troubled cells, i.e, cells on which the solution oscillates or behaves badly due to shocks or strong gradients have to be identified. Then, the polynomial of the troubled cell is reconstructed with the help of neighboring cells. This results in a local correction strategy.

In the following we will describe such strategies and how they can be applied in the sLdG context. Section \ref{sec:indicators} considers the indicator, i.e.~how do we decide which cells to modify, and section \ref{sec:modifier} then describes techniques for modification. As we will see, although these approaches work well, they are also relatively expensive on GPUs. Thus, in section \ref{sec:proposed-limiter} we will propose a new limiter that is computationally more efficient in the context of the sLdG scheme.

\subsection{Indicators \label{sec:indicators}}

We consider two methods here. The first strategy is a minmod approach which uses an idea based on finite volume slope limiters, see~\cite{Zhong2013} and the references therein.  The second is a scheme based on the extension of the polynomial in a cell to its immediate neighbors combined with a mean error indicator \cite{Fu2017}. 

We have chosen these two indicators from literature as they represent two classes of indicators. The minmod indicator is known to mark many cells as troubled and has thus the potential to introduce a significant amount of additional numerical diffusion. The mean error indicator based approach is known to be among the best and marks much fewer cells as troubled. It is therefore capable of preserving many more details of the density function under consideration. This also allows us to compare the two approaches and see how much effect additional numerical diffusion has in our problem. We refer to~\cite{Kochi2023,Qiu2005} for an overview and comparison of indicators from the literature. 

We will use the following notation. The polynomial in the cell that is to be checked is denoted by $u_0$. Further, $u_{-1}$ is the polynomial one cell to the left of $u_0$ and $u_{1}$ is the polynomial on cell to the right of $u_0$. We scale and shift these polynomials to the standard cells $I_{-1} = [-\frac{3}{2},-\frac{1}{2}]$, $I_0=[-\frac{1}{2},\frac{1}{2}]$ and $I_1=[\frac{1}{2},\frac{3}{2}]$, respectively. With a slight abuse of notation, we keep the same names, i.e.~$u_{-1}$, $u_{0}$ and $u_{1}$, for these polynomials.

\subsubsection*{Minmod}

In~\cite{Zhong2013}, the following minmod indicator is used. First, the mean 
\begin{equation}
\label{eq:compute_mean_in_cell}
\bar{u}_i = \int_{I_i} u_i\,dx, \quad i=-1,0,1,
\end{equation}
is computed. Next
\begin{align*}
u_0^+ &= u_0(\tfrac{1}{2})-\bar{u}_0, \\
u_0^- &= \bar{u}_0 - u_0(-\tfrac{1}{2}), \\
\Delta^+u &= \bar{u}_1-\bar{u}_0, \\
\Delta^-u &= \bar{u}_0-\bar{u}_{-1}, 
\end{align*}
are determined. Then, the minmod function 
 \[
     \text{minmod}(a,b,c) = \left\{
   \begin{aligned}
    &s\cdot\min(|a|,|b|,|c|) &\text{if } s = \sign(a)=\sign(b)=\sign(c) \\
    &0,  &\text{otherwise} \hspace{3.9cm} \\ 
   \end{aligned}
    \right.
 \]
 is applied. That is,
\begin{align*}
    \tilde{u}^+ &= \text{minmod}(u_0^+,\Delta^+u,\Delta^-u), \\
    \tilde{u}^- &= \text{minmod}(u_0^-,\Delta^+u,\Delta^-u),
\end{align*}
are computed. A cell is then marked as troubled whenever $\tilde{u}$ is different from the first argument of the minmod function.
Note that this method is completely parameter free.
It is well known that this limiter marks many cells as troubled and thus introduces a significant amount of additional numerical diffusion. In order to avoid this, the TVB modified minmod function~\cite{Cockburn1989} can be used. However, in this article we will use the original minmod function in order to compare methods that mark many cells as troubled compared to methods that mark few cells as troubled.

\subsubsection*{Mean error indicator}

We will now consider the indicator proposed in \cite{Fu2017}. The main idea of this approach is that the neighboring polynomials $u_{-1}$ and $u_{1}$ are naturally extended to the interval $I_0$. Then, the mean of the extended polynomials over $I_0$ is computed, i.e.,
\[
 \bar{\bar{u}}_i = \int_{I_0}u_i \,dx, \quad i=-1,0,1.
\]
Next, the average values $\bar{u}_{-1}$, $\bar{u}_0$ and $\bar{u}_1$ as in~\eqref{eq:compute_mean_in_cell} are determined. Then,
\begin{equation}
\label{eq:mean_err_indicator_orig}
 I_\Delta = \frac{|\bar{\bar{u}}_1-\bar{\bar{u}}_0|+|\bar{\bar{u}}_2-\bar{\bar{u}}_0|}{\max_i |\bar{u}_i|}
\end{equation}
is computed. A cell is then marked as troubled if $I_\Delta$ is bigger than a certain threshold (chosen as $0.5$ in the numerical experiments), which only depends on the degree of the polynomials. This method has the advantage that it can be extended to arbitrary dimensions and is one of the best indicators known in the literature, see~\cite{Fu2017,Kochi2023,Qiu2005}.

\subsection{Modifier \label{sec:modifier}}

Once a troubled cell is identified, it can be modified based on the polynomials from neighboring cells. In the following we will consider the modifiers in \cite{Zhong2013,Zhu2020,Zhu2016}, which are all mass conservative. We note note that we did not achieve satisfactory results using the modifier from~\cite{Zhu2016} and thus omit it here.

\subsubsection*{Simple WENO}

The method introduced in~\cite{Zhong2013} is a WENO scheme that proceeds as follows. First, the smoothness coefficients
\begin{equation}
\label{eq:smoothness_indicators}
\beta_l = \sum_{s=1}^k\int_{I_l}\left(\frac{\partial^s}{\partial x^s} u_l(x)\right)^2 \,dx,
\end{equation}
are computed, see \cite{Jiang1996}). Then, the nonlinear weights 
\[
 \bar{\omega}_l = \frac{\gamma_l}{(\varepsilon + \beta_l)^2},
\]
where the linear weights $\gamma_l$ are defined as 
$\gamma_{-1}=0.001$, $\gamma_{0}=0.998$, $\gamma_{1}=0.001$ and $\varepsilon=10^{-6}$ is used in order to avoid problems with division by zero.
Next, these weights are normalized
\begin{equation} \label{eq:omegal}
 \omega_l = \frac{\bar{\omega}_l}{\sum_s \bar{\omega}_s}
\end{equation}
such that we have $\omega_{-1} + \omega_0 + \omega_{1}=1$.
Finally, the modified polynomial is constructed as follows
\[
 u_0^{\text{new}}(x) = \omega_{-1} \left(u_{-1}(x) - \bar{u}_{-1} +\bar{u}_0\right) + \omega_0 u_0(x) + \omega_{1}\left(u_{1}(x) - \bar{u}_{1} + \bar{u}_0\right).
\]

\subsubsection*{Line WENO}

Next, we describe the limiter developed in~\cite{Zhu2020}, which in order to distinguish it from the limiter in the previous section we will call line WENO in this article.
Here, again the mean values of the polynomials~\eqref{eq:compute_mean_in_cell} are computed. Next, we find the unique linear polynomials $p_{-1}$ and $p_1$ that satisfy
\begin{align*}
 &\int_{I_s}p_{-1}(x)\,dx = \bar{u}_s, \quad  s=-1,0 \\
 &\int_{I_s}p_1(x)\,dx = \bar{u}_s, \quad  s=0,1. 
\end{align*}
Then, with the help of three linear weights $\gamma_{-1}$, $\gamma_0$ and $\gamma_1$, the polynomial 
\[
p_0(x) = \frac{1}{\gamma_0}u_0(x) - \frac{\gamma_{-1}}{\gamma_0}p_{-1}(x) - \frac{\gamma_1}{\gamma_0}p_1(x)
\]
is generated. These linear weights satisfy $\gamma_{-1} + \gamma_0 + \gamma_1 = 1$. In the numerical experiments we have chosen $\gamma_1 = \gamma_{-1} = 0.45$ and $\gamma_0=0.1$. In a similar fashion as for the simple WENO modifier, the nonlinear weights $\omega_l$ are computed as follows
\[
 \tau = \left(\frac{|\beta_0-\beta_{-1}| + |\beta_0-\beta_1|}{2}\right)^2,
\]
where $\beta_l$ are defined in~\eqref{eq:smoothness_indicators}. Then, the nonlinear weights are given as
\[
 \bar{\omega}_l = \gamma_l\left(1+\frac{\tau}{\varepsilon + \beta_l}\right), \quad l=-1,0,1.
\]
We again use \eqref{eq:omegal} to obtain $\omega_l$.
Then, the modified polynomial is constructed as follows
\[
 u_0^{\text{new}}(x) = \omega_{-1}p_{-1}(x) + \omega_0p_0(x) + \omega_1p_1(x).
\]

\begin{remark}
\label{rmk:downside_classic_limiter}
Although these limiters are local and only require data from two neighboring cells of the troubled cell, they are still not ideal from a computational point of view in the context of the sLdG method.  First, these limiters are applied after (or before) the advection step is performed and thus require that the updated values are already available in all neighboring cells. If this is implemented as a postprocessing step, then additional reads and writes of all degrees of freedom from memory are necessary. Since the sLdG scheme is memory bound this can be as or even more costly than the advection step itself.

Second, in a distributed memory context (i.e. multiple GPUs or multiple multi-core CPUs that are connected via a network) where the problem is parallelized to several subdomains, boundary data has to be transferred. As shown in \cite{Einkemmer2022}, sLdG methods only require boundary data in the upwind direction. The described limiters, however, require boundary data from both directions. This is a massive drawback on large supercomputers, where the data transfer is usually the performance bottleneck \cite{Einkemmer2022}. 
\end{remark}

\subsection{Computationally efficient limiter for the sLdG scheme \label{sec:proposed-limiter}}

In order to overcome the downsides pointed out in Remark~\ref{rmk:downside_classic_limiter}, we propose a limiter which can be directly applied in the advection step. Also in this case we first check if a cell is troubled. If this is the case, the polynomial in the trouble cell is modified. However, in contrast to the limiter described sections above, only the data from two input cells and one output cell is used (see the form of the sLdG algorithm in equation~\eqref{eq:sldg_implementation}).

In the following, we denote with $u_l$ and $u_r$ the left and right input cell scaled and shifted to $I_l=[0,1]$ and $I_r=[1,2]$, respectively. In addition, we denote the projected polynomial on the interval $I_p=[\alpha,1+\alpha]$, $\alpha \in [0,1)$ by $u_p$.

\subsubsection*{Indicator} 

We adjust the mean error indicator to this new setting. We naturally extend $u_p$ to the intervals $I_l$ and $I_r$. Then we compute
\begin{align*}
\bar{u}_{p,l} = \int_{I_l} u_p(x)\,dx \quad \text{ and } \quad \bar{u}_{p,r} = \int_{I_r} u_p(x)\,dx
\end{align*}
After computing the mean of $\bar{u}_l$ and $\bar{u}_r$ as in~\eqref{eq:compute_mean_in_cell}, we can, similar to equation~\eqref{eq:mean_err_indicator_orig}, compute 
\begin{equation}
  I_\Delta^{\text{mod}} = \frac{|\bar{u}_{p,l}-\bar{u}_l| + |\bar{u}_{p,r}-\bar{u}_r| }{\max\{|\bar{u}_l|,|\bar{u}_r|\}}.
\end{equation}
A polynomial $u_p$ is then marked as troubled if $I_\Delta^{\text{mod}}$ is bigger than a given threshold (chosen as $0.5$ in the numerical experiments).

\begin{figure}
\centering 
\includegraphics[width=0.4\textwidth]{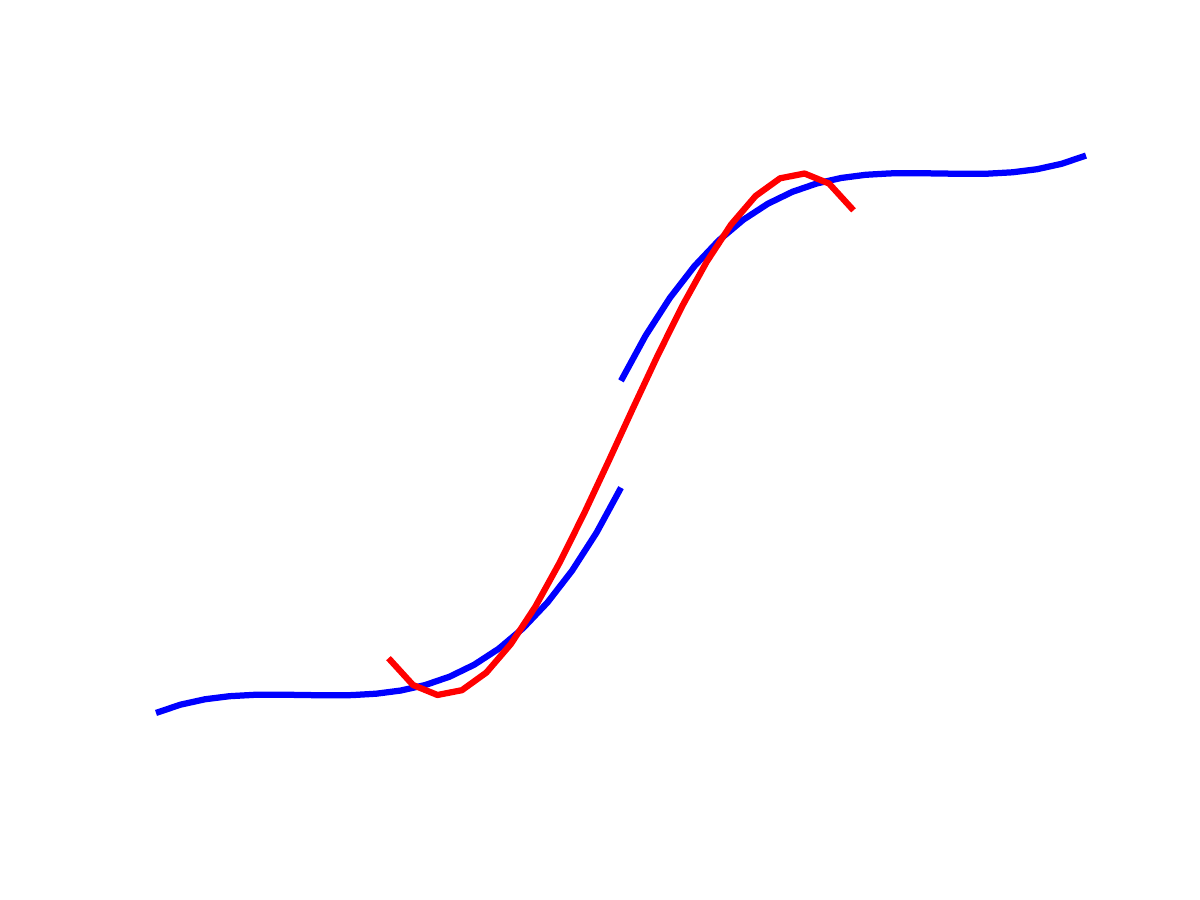}%
\includegraphics[width=0.4\textwidth]{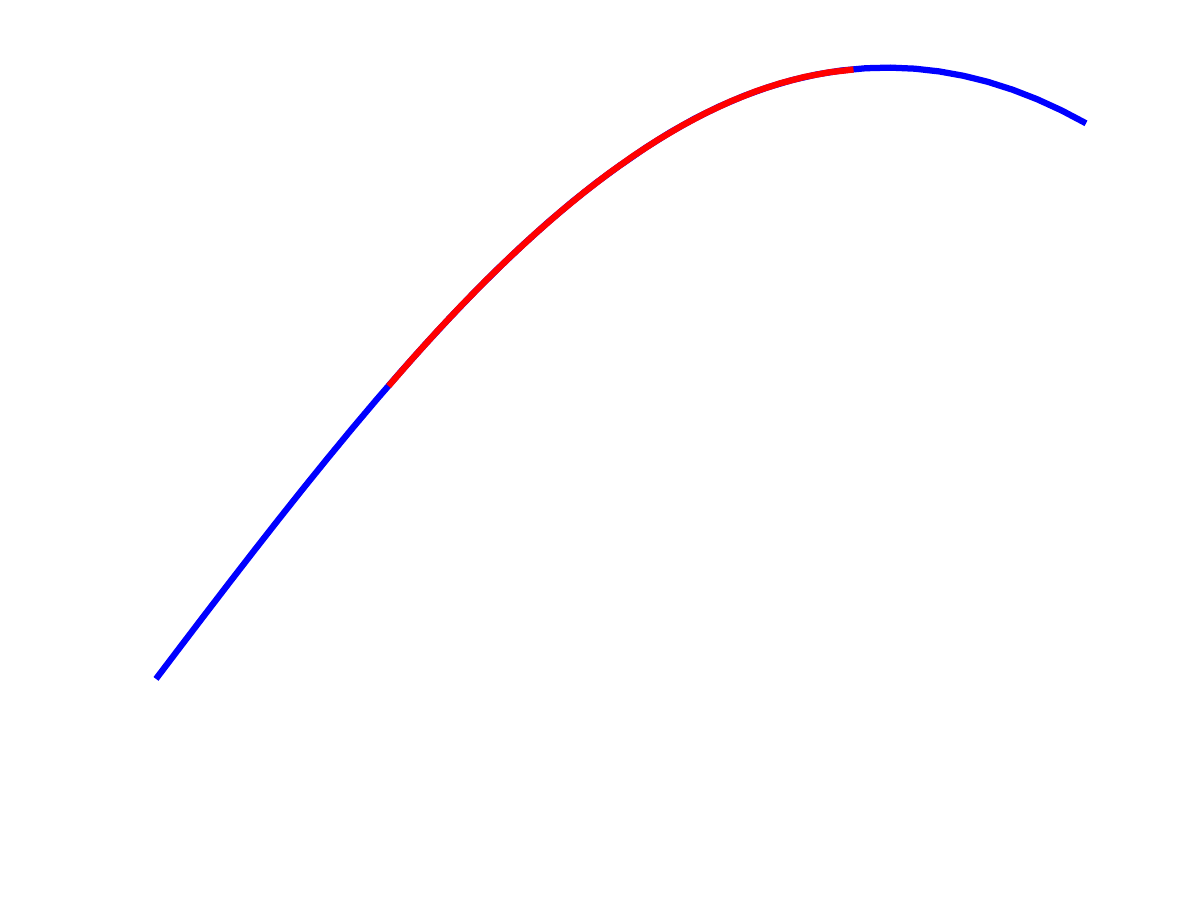}
\caption{Left: Projection applied in a region with a steep derivative. Right: Projection applied in a smooth region. The blue curves are the input polynomials ($u_l$ and $u_r$) and the red curves are the projections $u_p$. Here polynomials of degree $k=3$ are used.}
\label{fig:projection_demo}
\end{figure}

Before proceeding, let us remark that projections can introduce oscillations in regions with a steep gradient (or discontinuity) if not enough resolution is provided, see figure~\ref{fig:projection_demo}. The proposed troubled cell indicator works excellent in this case as the extension of the projection differs significantly from the input cells, while on smooth regions the extension follows the input cell accurately.

\subsubsection*{Modifier} 

In order to modify a troubled cell, we apply a maximum-minimum limiter on $u_p$ compared to the input $u_l$ and $u_r$. That is, we want to have
\begin{align*}
u^{\text{max}}_p = \max_{x\in I_p}\, u_p(x) &\leq \max \left( 
\max_{x\in[\alpha,1]}u_l(x),
\max_{x\in[1,1+\alpha]}u_r(x)\right) = u^\text{max}_{l,r}, \\ 
u^{\text{min}}_p = \min_{x\in I_p}\, u_p(x) &\geq \min \left( 
\min_{x\in[\alpha,1]}u_l(x),
\min_{x\in[1,1+\alpha]}u_r(x)\right) = u^\text{min}_{l,r}.
\end{align*}
In order to achieve this, we compute
\begin{equation}
 \theta = \min\left(\left|\frac{u^{\text{max}}_{l,r}-\bar{u}_p}{u^{\text{max}}_p-\bar{u}_p}\right|,\left|\frac{u^{\text{min}}_{l,r}-\bar{u}_p}{u^{\text{min}}_p-\bar{u}_p}\right|,1\right)
\end{equation}
and modify the output as follows
\[
 u_p^{\text{new}}(x) = \theta\left(u_p(x) - \bar{u}_p\right) + \bar{u}_p
\]
This modifier can be seen as an extension of global maxima and minima limiter or positivity preserving limiters, see~\cite{Rossmanith2011,Qiu2011}. 
Let us remark that although it seems natural to apply this modifier to every cell (with the added benefit that the solution remains positive), numerical tests show that too much artificial numerical diffusion is introduced in this case. 

\begin{remark}
 For SOL simulations, one primary quantity of interest is the flux at $x=L$, which is defined as 
\begin{equation}
\label{eq:wrong_flux}
 j(t) = \int_\mathbb{R}vf(t,L,v)\,dv.
\end{equation}
The limiters used in this article have a negative impact on this quantity as they have a diffusive behavior on the density function close to the zero inflow boundary.
Therefore, we compute the flux by taking the integral only over the outflow velocities, i.e.
\begin{equation}
\label{eq:flux} 
 j(t) = \int_{\mathbb{R} \geq 0}vf(t,L,v)\,dv.
\end{equation}
In the continuous setting both integrals are exactly equivalent. However, in the actual numerical implementation the effect of the limiter on the inflow boundary is avoided by using equation \eqref{eq:flux}.
\end{remark}

\section{Adaptive mesh refinement \label{sec:amr}}

The scrape of layer (SOL) simulations considered in this paper have two distinct feature. First, particles with high velocity (particularly electrons) leave the computational domain relatively early in the simulation. Thus, if a fixed velocity domain is chosen a significant amount of computational power is wasted doing no useful computation. This is true, in particular, since the main heat flux at the divertor is determined by the bulk ions (at which points the tail of the original density function has already been lost). Second, close to the boundary a spatial boundary layer, the so-called plasma sheath develops as the faster electrons result in an excess of negative charge close to the wall. This region is characterized by sharp gradients that need to be resolved numerically. Both of these phenomenon can be treated by adaptive mesh refinement. The approach chosen in this work is described in the remainder of this section.

\subsection{Adaptive velocity domain adjustment \label{sec:adaptivevel}}

Particles initially have a wide range of velocities. Some particles travel with high speed along the magnetic field lines and collide with the divertor plates of the tokamak, which we assume here to be perfectly absorbing. These particles are thus lost. The result is that the support of the density function $f_{\alpha}$ in the velocity direction shrinks as time progresses. Therefore, after some time, most of the degrees of freedom just store zeros. In order to avoid this and increase the effective resolution when the bulk of the density function hits the divertor, we adjust the domain in the velocity direction adaptively. 

While changing from one grid to another one is relatively cheap within the discontinuous Galerkin framework, we still want to avoid doing this too often. Therefore, at the end of each time step, we check the following condition
\begin{equation}
\label{eq:v_adjustment}
\begin{aligned}
f_\alpha(t,x,v=v_\text{max}(1-p-\gamma)) &< \text{tol} \\
f_\alpha(t,x,v=v_\text{min}(1-p-\gamma)) &< \text{tol},
\end{aligned}
\end{equation}
where the parameter $p$ defines by how much the velocity domain bounds are to be reduced (in the numerical simulation we have chosen $p=0.05$, which corresponds in a reduction of 10\% of the numerical support). The parameter $\gamma$ is a safety factor (in the numerical simulations we have chosen $\gamma=0.05$). Verifying this condition is a relatively cheap operation as only two 1d stripes of $f_\alpha$ which lie continuous in memory have to be considered. 

If the criteria in~\eqref{eq:v_adjustment} are satisfied (we have set $\text{tol}=10^{-14}$), $f_\alpha$ is projected onto the smaller domain
\[
 [-L,L]\times[-v_\text{max}(1-p),v_\text{max}(1-p)].
\]
Within the semi-Lagrangian framework this projection can be done simply by locally projecting the piecewise polynomial to the new cells. Similar as for the advection step in ~\eqref{eq:sldg_implementation} this can be cast in the form of small matrix-vector products of adjacent 1d cells. These matrices are computed for every velocity cell on the fly and are then used for each degree of freedom in the $x$ direction. As a consequence and because we do not need to this in every time step, no significant computational overhead is introduced by performing this adaptive domain adjustment in the velocity direction.

This is highly beneficial for SOL layer simulations as initially a large velocity domain is available to capture the fast particles, while later when the bulk of the plasma hits the wall the effective resolution available compared to a non-adaptive velocity domain is drastically increased. We will study this in more detail in section \ref{sec:numerical}.

\subsection{Refinement close to the wall \label{sec:refwall}}

In order to resolve the step gradient and possible turbulence close to the wall, a high resolution is required. However, running the entire simulation with such a resolution would result in a prohibitive computational cost. Since the size of the sheath is relatively static, we simply split the domain into three blocks in the $x$-direction. That is, we have two blocks close to the wall with cell size $\Delta x_\text{fine}$ and one interior block where the cell size is $\Delta x_\text{coarse}$. For simplicity, we define the relation $\Delta x_\text{coarse} = m \Delta x_\text{fine}$, where $m\in\mathbb{N}$.

Note that the transport in the velocity direction is not affected by this choice. However, the implementation of the transport in the spatial direction would have to account for the change in cell size. In particular, when transitioning from the coarse to the fine grid more than two adjacent cells would be involved in an advection step. We avoid this by using ghost cells for which the advection can be computed according to equation \eqref{eq:sldg_implementation} and using the fast implementation that is already available. As a postprocessing step, the transfer from the fine to the coarse mesh and vice versa is performed. To transfer data from the fine to the coarse grid, first $m\cdot\text{CFL}_{\text{coarse}}$ polynomials on the fine grid are projected to $\text{CFL}_\text{coarse}$ cells of the coarse grid, where $\text{CFL}_\text{coarse}$ is the CFL number in the interior block. Once this projection has been performed, only CFL\textsubscript{coarse} boundary cells have to be transferred. This is beneficial because it reduces the amount of data to transfer, which is important in high-performance computing. 
On the other hand, going from the coarse to the fine grid is straightforward polynomial evaluation. In this case also $\text{CFL}_\text{coarse}$ boundary cells need to be transferred. The outlined approach also has the advantage that it allows us to run the different blocks e.g.~on different GPUs or different nodes of a supercomputer.

\section{Numerical experiments \label{sec:numerical}}

In this section we consider two 1x1v simulations of the SOL. In section \ref{sec:blob} we consider a blob initial value (i.e.~the plasma localized in the interior of the domain). Such a setting is often used as a simplified model for plasma that is expelled from the core of the device in an ELM (edge-localized mode) burst. In section \ref{sec:injection} we then consider a setting where the plasma is continuously injected into the SOL. Mathematically this is modeled by an appropriate source term. Finally, we investigate the performance of the implementation of the various limiters in section \ref{sec:performance}.

We consider a domain of length $L=200$, a time step size $\Delta t=0.1$, and the sLdG scheme of order $4$. The mass ratio is set to $\mu = 400$. This is similar to the configuration chosen in~\cite{Bourne2023} and fairly typical what is reported in the literature. We note, however, that the physical mass ratio for an electron-hydrogen plasma is $\mu \approx 1836$ and for a realistic tokamak $L \approx 10^6$. However, smaller values for $L$ can be used without too much loss in accuracy, see~\cite{Coulette2016} where $L = 10^3$ is chosen. All simulations have been conducted using our open source SLDG code\footnote{\url{https://bitbucket.org/leinkemmer/sldg}}.

\subsection{Initial blob \label{sec:blob}}

\begin{figure}
 \includegraphics[width=\linewidth]{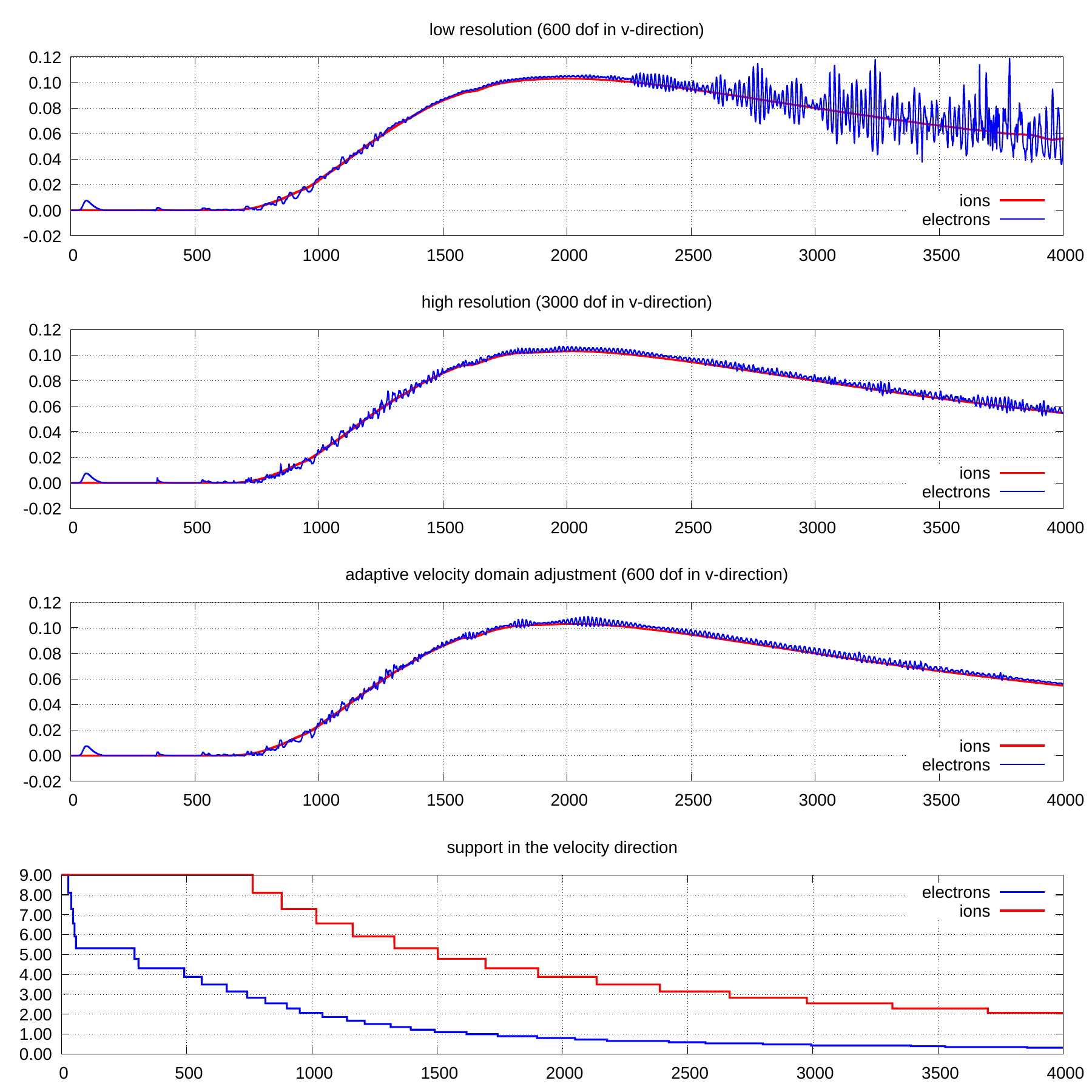}
 \caption{Particle flux at the right boundary~\eqref{eq:flux} with coarse velocity resolution (top), fine velocity resolution (second from top), and coarse resolution with the adaptive velocity adjustment algorithm (second from bottom). In addition, the bottom plot shows how the largest speed that can be found in the system decreases in time. In all simulations $1200$ degrees of freedom are used in the spatial direction.}
\label{fig:vdomain+problems} 
\end{figure}

We consider the following initial condition
\[
f_\alpha(t=0,x,v) = \frac{1}{\sqrt{2\pi}}\exp\left(-\frac{x^2}{2\sigma^2}\right)\exp\left(-\frac{v^2}{2}\right),
\]
where $\sigma=0.1L$. In figure~\ref{fig:vdomain+problems} we report the particle flux~\eqref{eq:flux} at the right boundary. Oscillations are primarily an issue for the electron flux. The reason being that electrons are lighter and thus move with higher velocity, which results in a more filamented phase space. Consequently, the amount of oscillations can be significantly reduced by increasing the velocity resolution (see second from top plot). Clearly this comes at greatly increased computational cost, which we want to avoid. We also observe (see top of Figure~\ref{fig:vdomain+problems}) that the oscillations increase with time. This points to the fact that as the fast electrons are lost the effective velocity space is covered by only a couple of cells. The numerical solution is then severely underresolved which leads to oscillations. As we can see from the second from bottom plot, the adaptive velocity adjustment described in section \ref{sec:adaptivevel} allows us to simulate the problem with much less degrees of freedom. In the bottom plot of figure~\ref{fig:vdomain+problems} we report the time evolution of the numerical support in the velocity domain. For the electrons the velocity domain is reduced by approximately a factor of $20$ from the start to the end of the simulation.

\begin{figure}
 \includegraphics[width=\linewidth]{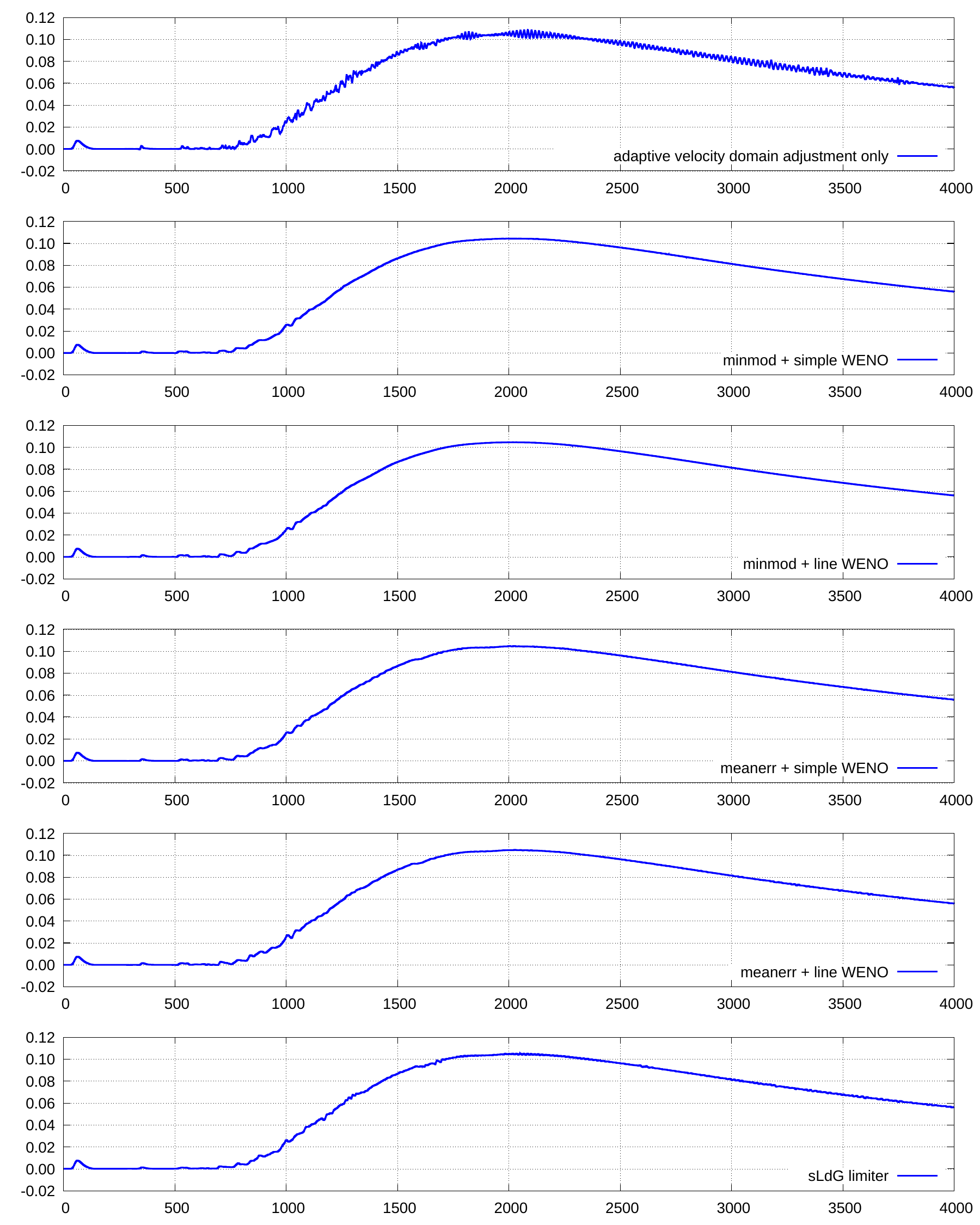}
 \caption{Electron flux for the blob initial data at the right boundary~\eqref{eq:flux} with 1200 degrees of freedom in the spatial direction, 600 degrees of freedom in the velocity direction, and using the adaptive velocity adjustment. The different limiters used are indicated in the legends of the plots. \label{fig:limiter}}
\end{figure}

We also observe that for the resolutions studied here the ion flux is smooth.  However, even with relatively high resolution or with the adaptive velocity adjustment enabled, some oscillations in the electron flux remain. Hence, from now on, we will be mainly concerned with (the more challenging) electron flux. In order to suppress these oscillations, we thus turn our attention to the limiters described in section \ref{sec:limiter} . The results for different indicators and modifiers and the newly proposed sLdG limiter are shown in Figure \ref{fig:limiter}. We observe that all limiters successfully remove the majority of the oscillations in the electron flux, while not having any negative impact on the physical features of the solution.

\begin{figure}
\centering
 \includegraphics[width=0.49\linewidth]{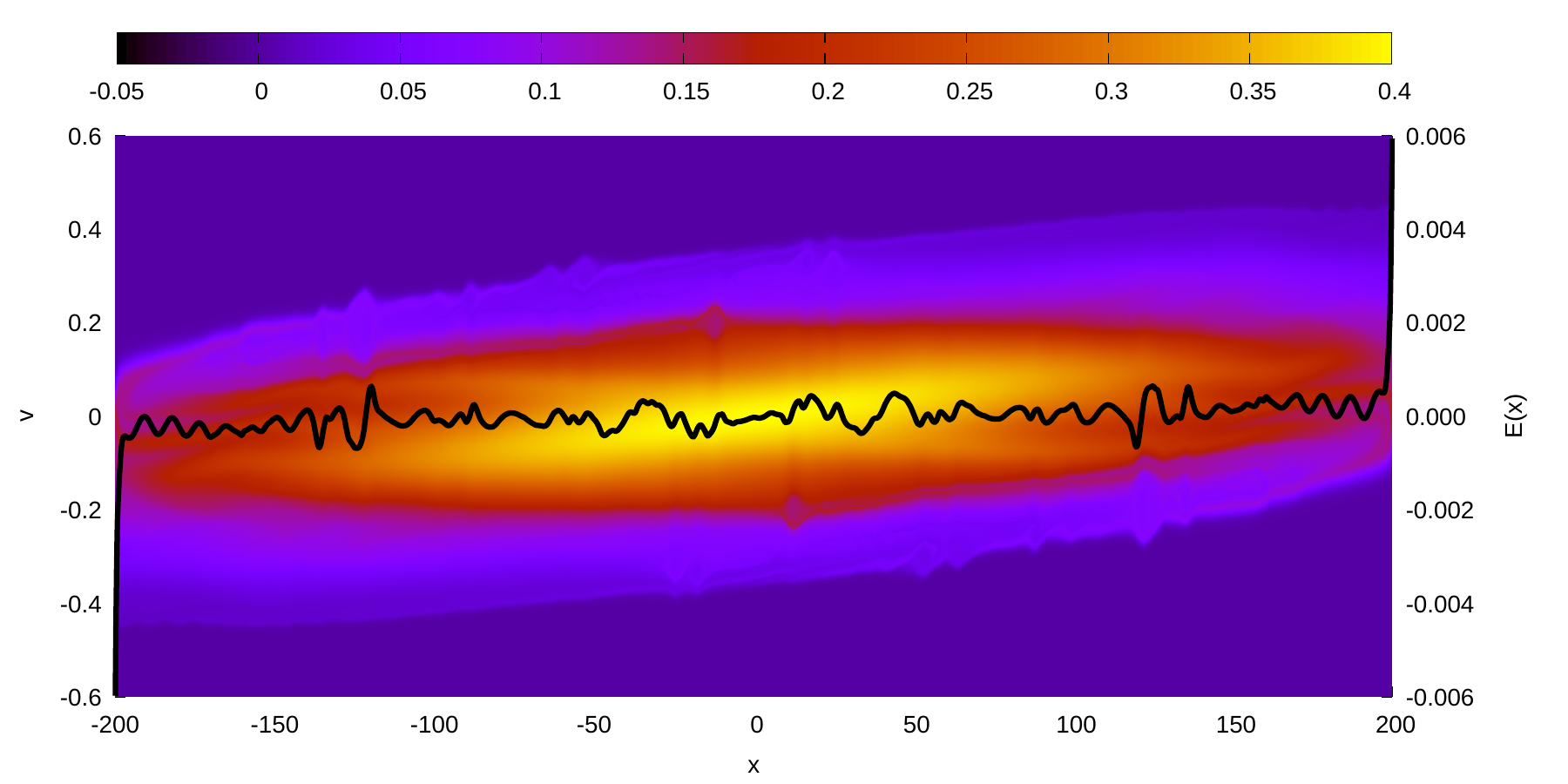}
 \includegraphics[width=0.49\linewidth]{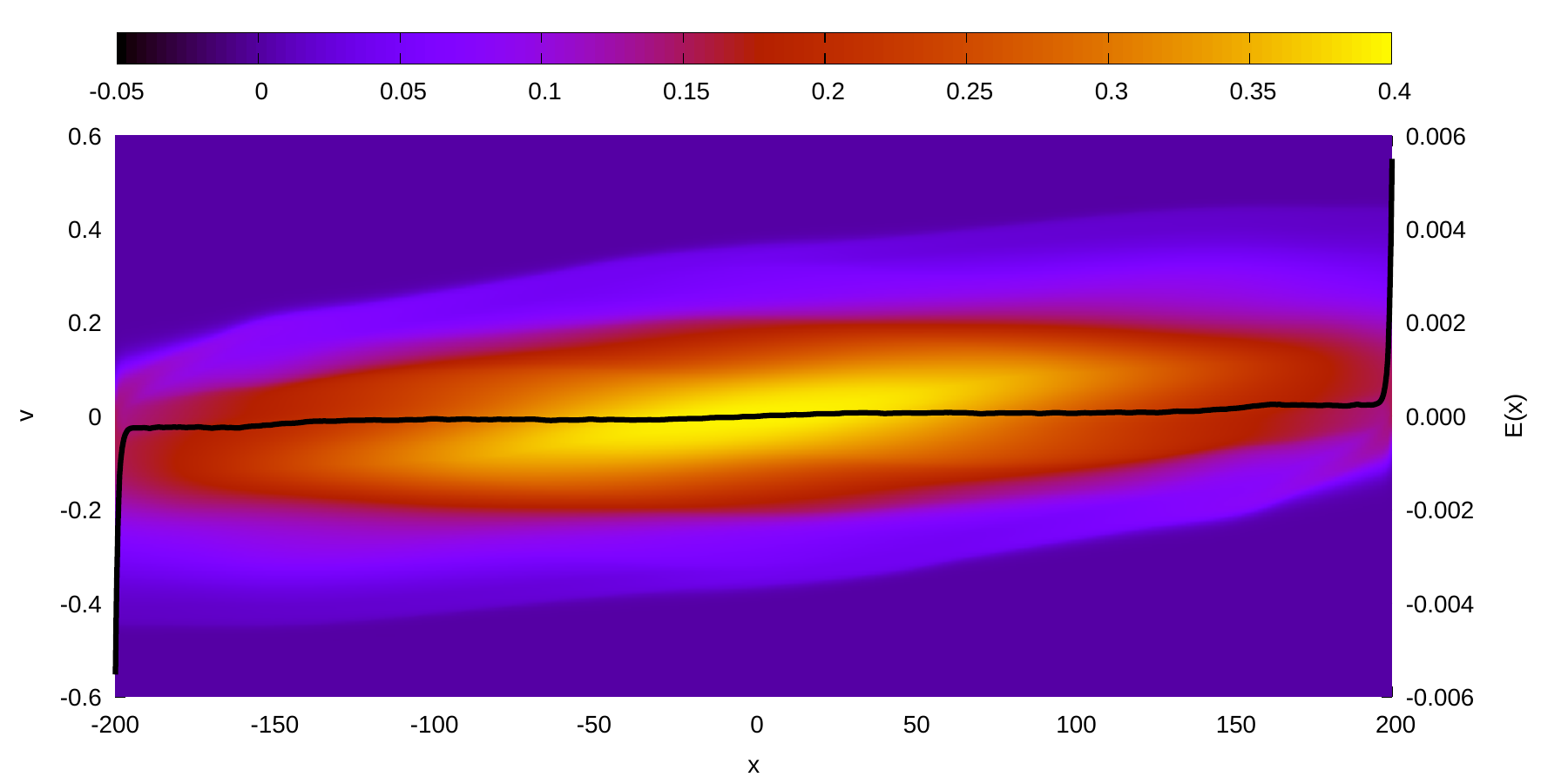} \\
 \includegraphics[width=0.49\linewidth]{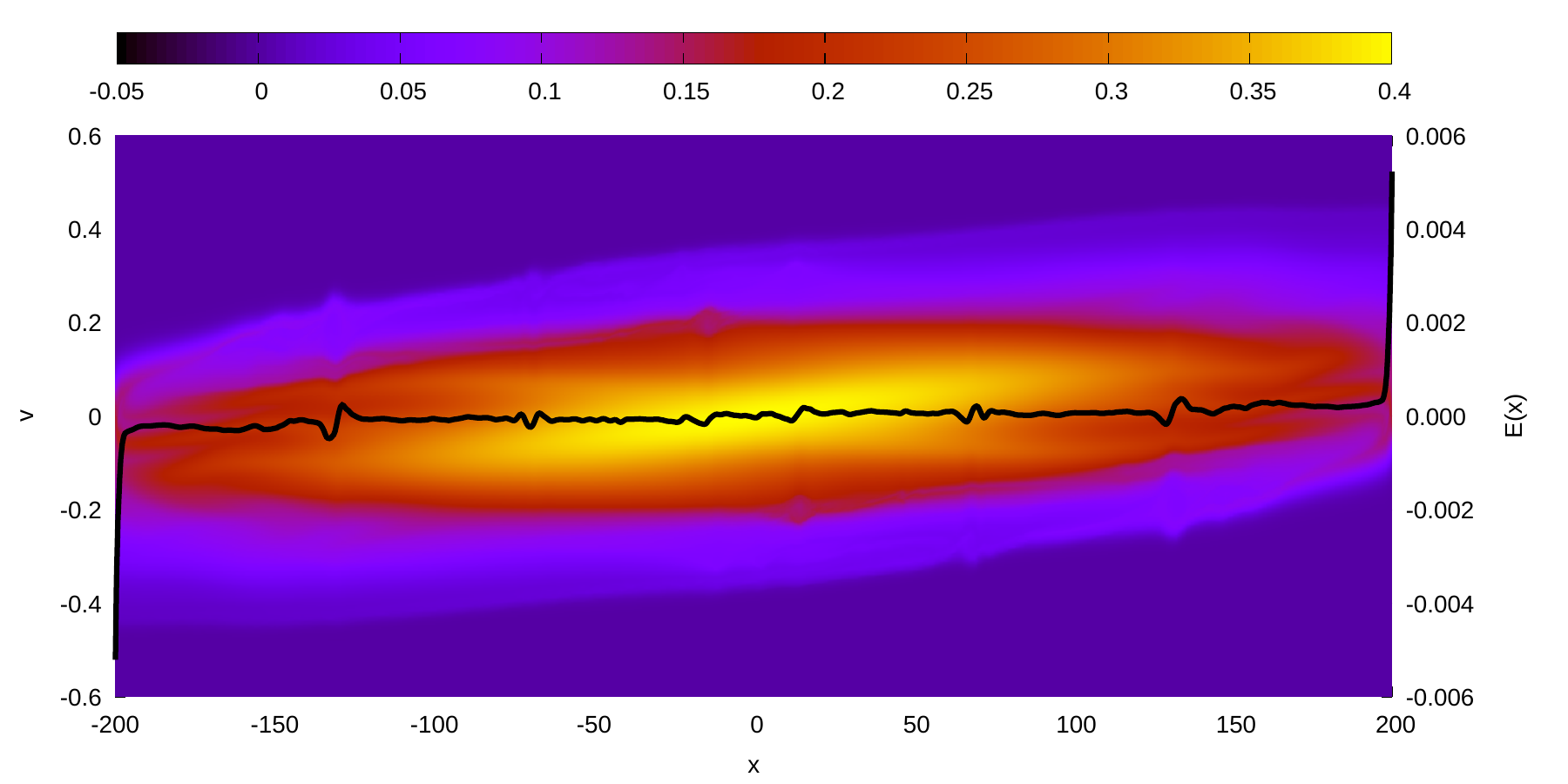}
 \includegraphics[width=0.49\linewidth]{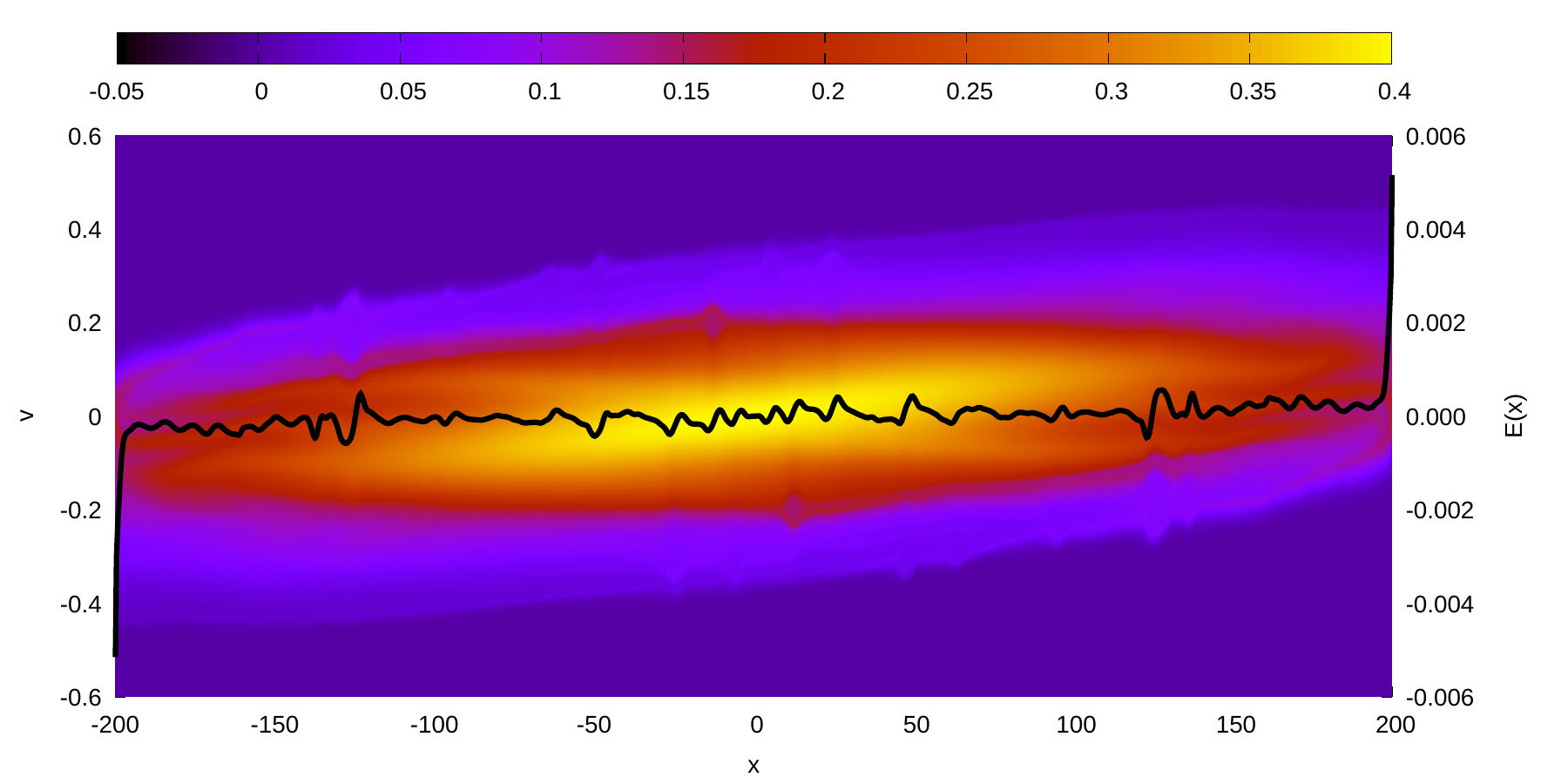}
 \caption{The electron particle density for the blob initial data at time $t=2000$ is shown for different limiters. 
Top left: reference solution, meanerr indicator and line WENO modifier, $2400 \times 1200$ degrees of freedom. Top right: minmod indicator and line WENO modifier, $1200 \times 600$ degrees of freedom. Bottom left: meanerr indicator and line WENO modifier, $1200 \times 600$ degrees of freedom. Bottom right: the newly proposed sLdG limiter, $1200 \times 600$ degrees of freedom. In all simulations the adaptive velocity adjustment is used.}
 \label{fig:density_ii}
\end{figure}

In Figure \ref{fig:density_ii} we show the density function and the electric field for selected limiters. Here we see significant differences. In particular, the minmod indicator is very diffusive and results in a perfectly smooth field and density function. Limiters that are based on variants of the meanerr indicator are significantly less diffusive, but correspondingly show more oscillations. This is also true for the proposed sLdG limiter. The main reason for this is that fewer cells are marked as troubled. We note, however, that for the reference solution, which is computed with increased resolution and limiter enabled, some oscillations do persist. Thus, the oscillations observed are not purely a numerical artifact, but a physical feature of the solution. In fact, the sLdG limiter matches the reference solution best, followed by the meanerr indicator, and then the minmod indicator. We also remark that in all cases oscillations at the sharp gradient that occur in the plasma sheath are avoided. Consequently, all limiters give relatively similar results close to the wall as we have seen in figure~\ref{fig:limiter}.

\begin{figure}
\centering 
\includegraphics[width=0.9\linewidth]{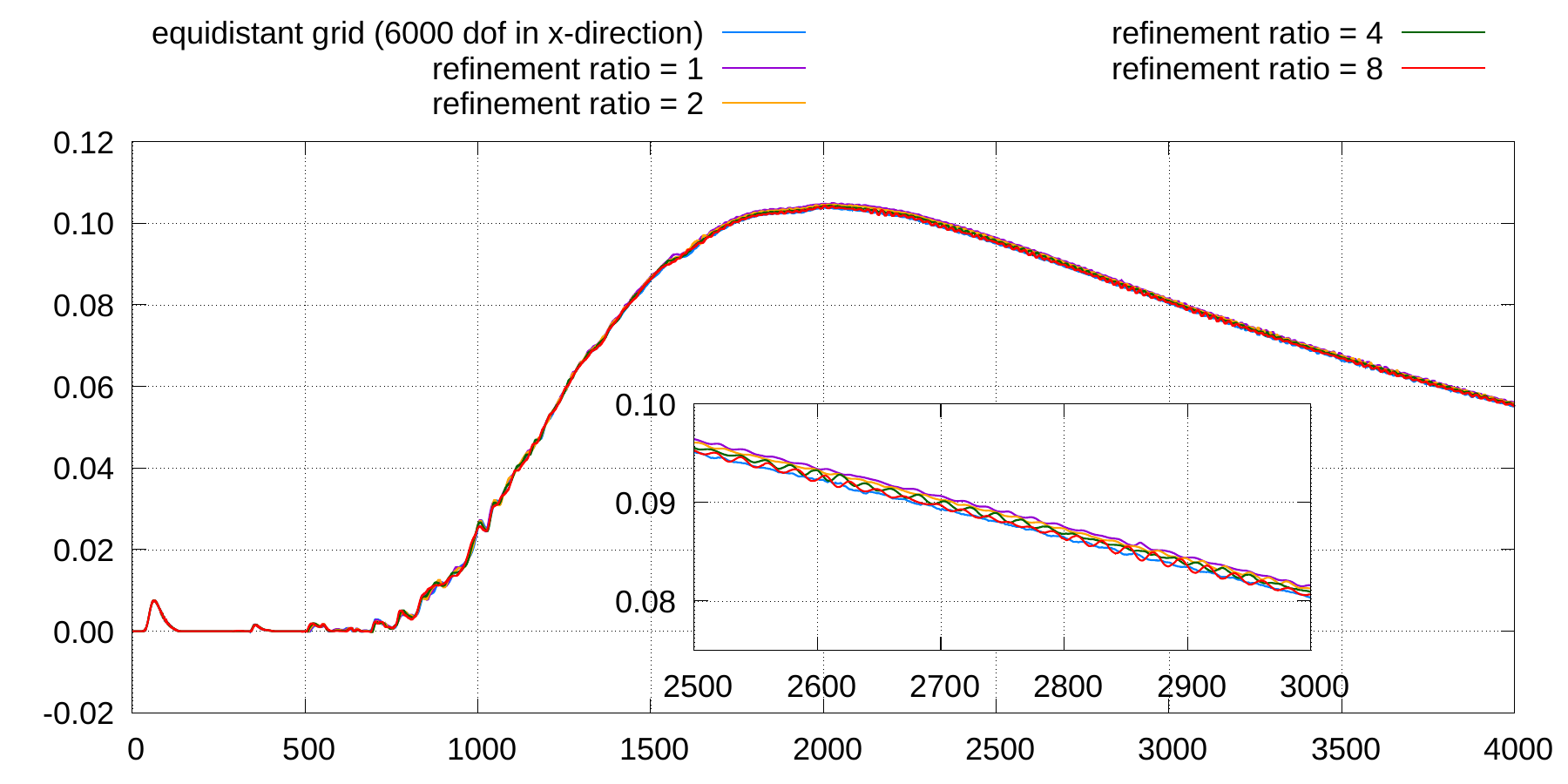}
\caption{Electron particle flux for the blob initial data using different refinement ratios close to the wall. The meanerr indicator and the line WENO modifier are applied and the total (i.e.~in all blocks) number of degrees of freedom is $1200 \times 600$. As a comparison a fine x-grid (with $6000 \times 600$ degrees of freedom) simulation with the same limiter is shown.}
\label{fig:xdom-adjustment}
\end{figure}

Before proceeding, let us turn our attention to the mesh refinement close to the wall. As outlined in section \ref{sec:refwall} we split our computational domain into three block. The refinement ratio ($m$ in section \ref{sec:refwall}) gives the ratio between the cells size for the coarse grid that covers the middle of the domain to the fine grid that covers the two areas close to the boundary. Here we have chosen the same number of degrees of freedom for each block. That is, as the refinement ratio increases the physical size of the two boundary blocks decreases, while the total number of degrees of freedom is held constant. In figure \ref{fig:xdom-adjustment} we compare these results with a reference solution that uses $6000$ degrees of freedom in the spatial direction. We observe that by increasing the refinement ratio the results of the simulations converge to the simulation with the fine resolution.

\subsection{Continuous injection of plasma \label{sec:injection}}

In this case we assume that the SOL is initially in a vacuum state. That is, we set $f_{\alpha}(0,x,v) = 0$. A localized blob of plasma then is injected using the source term
\[
 S(t,x) = \frac{H(t_0-t)}{\sqrt{2\pi} \cdot t_0}\exp\left(-\frac{x^2}{2\sigma^2}\right)\exp\left(-\frac{v^2}{2}\right),
\]
where $H$ is the Heaviside function and $t_0=2000$.

\begin{figure}
\centering
 \includegraphics[width=\linewidth]{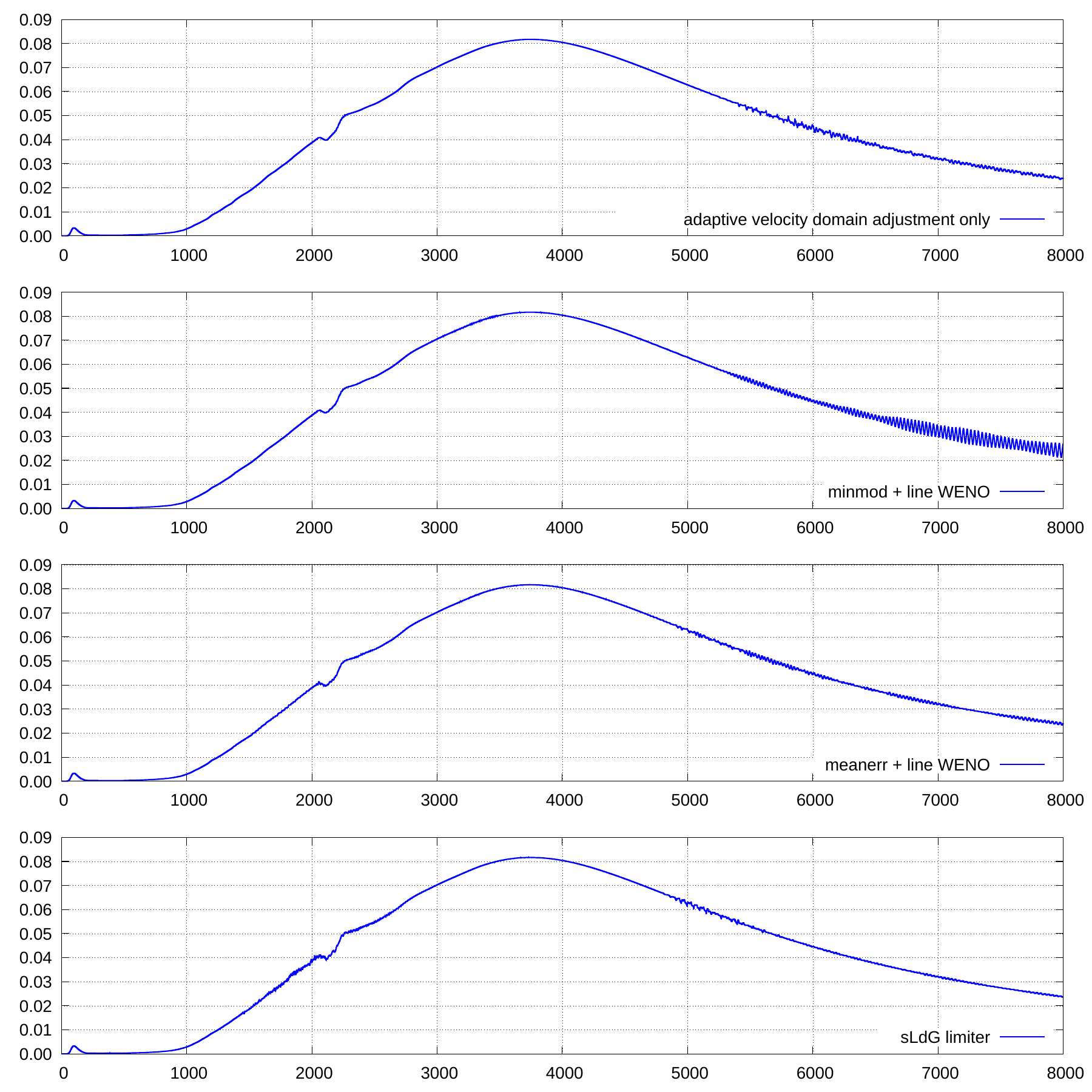}
 \caption{Electron flux for continuously injected plasma at the right boundary~\eqref{eq:flux} with $1200 \times 600$ degrees of freedom, using the adaptive velocity adjustment, and a refinement ratio of $8$ close to the wall. The different limiters used are indicated in the legends of the plots.  \label{fig:flux_cs}}
\end{figure}

In Figure \ref{fig:flux_cs} we show the electron flux at the right boundary for different limiters. Interestingly, using the minmod indicator substantially increases the amplitude of the oscillations observed in the numerical solution for times $t>5500$. Note that since the source term is active until $t=2000$, we have to run the simulation for longer ($t=8000$ compared to $t=4000$, as we have used before). Both the meanerr based limiters and the newly proposed sLdG limiters only show small amplitude oscillations and thus perform satisfactorily. For long times the newly proposed sLdG limiters suppresses the oscillations better than the other limiters, but for short times the meanerr based limiters are superior. 

\begin{figure}
\centering
 \includegraphics[width=0.49\linewidth]{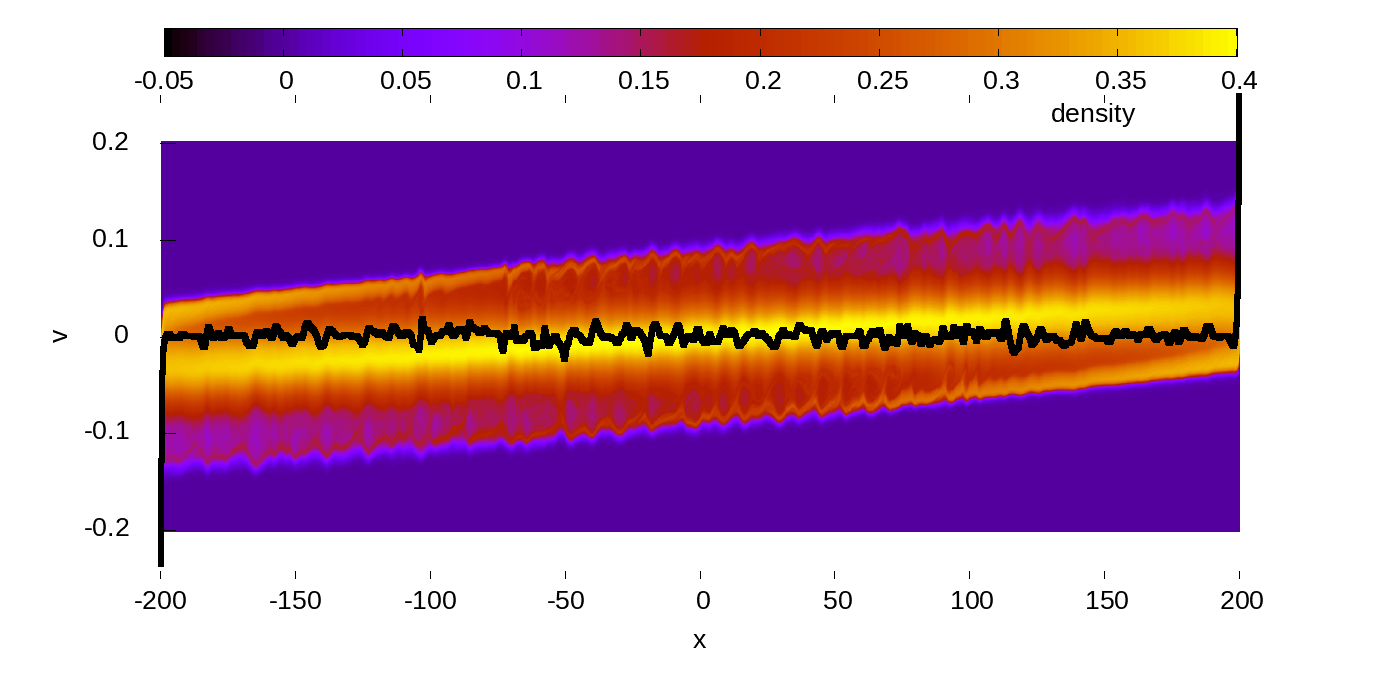}
 \includegraphics[width=0.49\linewidth]{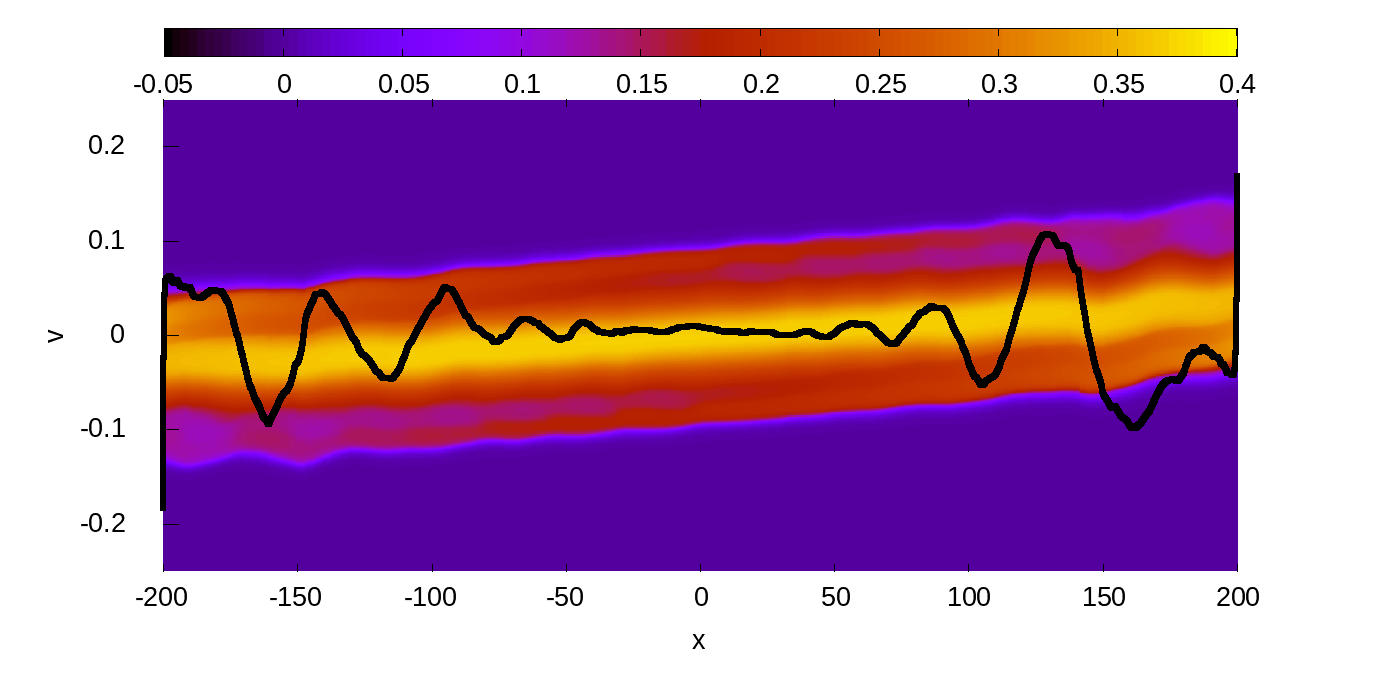}
 \includegraphics[width=0.49\linewidth]{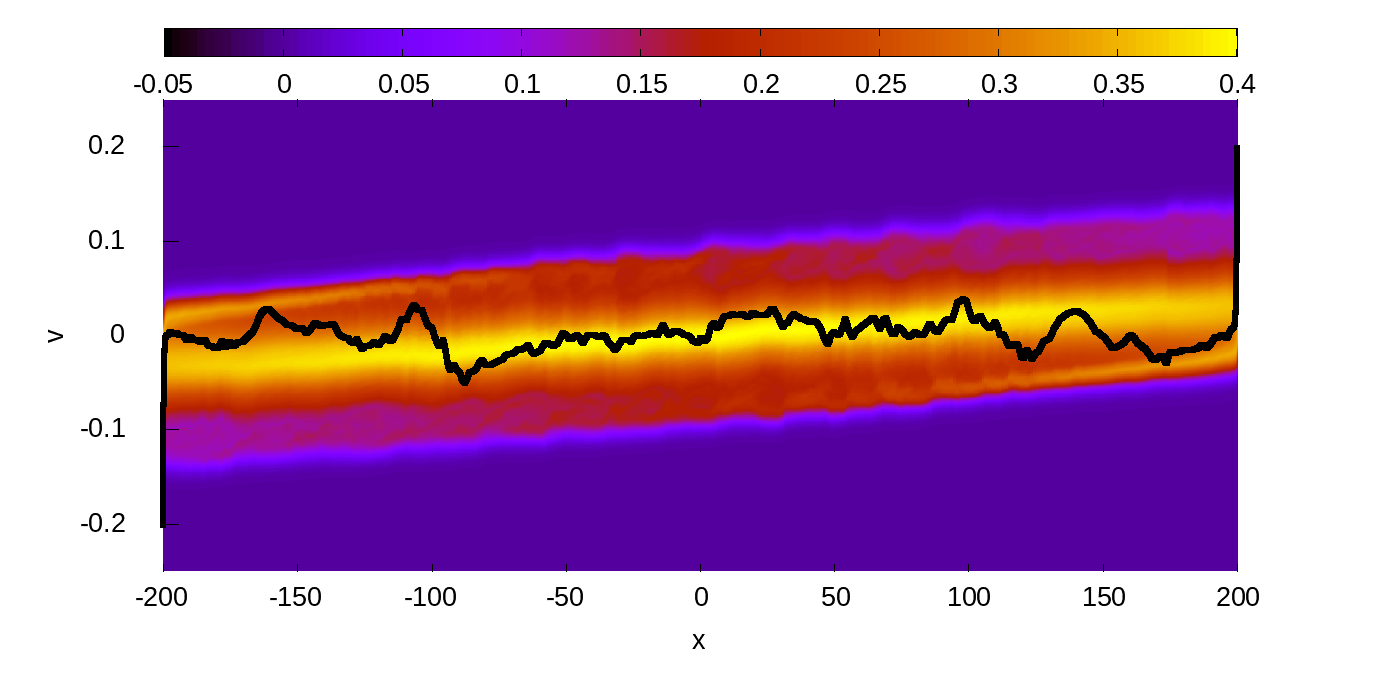}
 \includegraphics[width=0.49\linewidth]{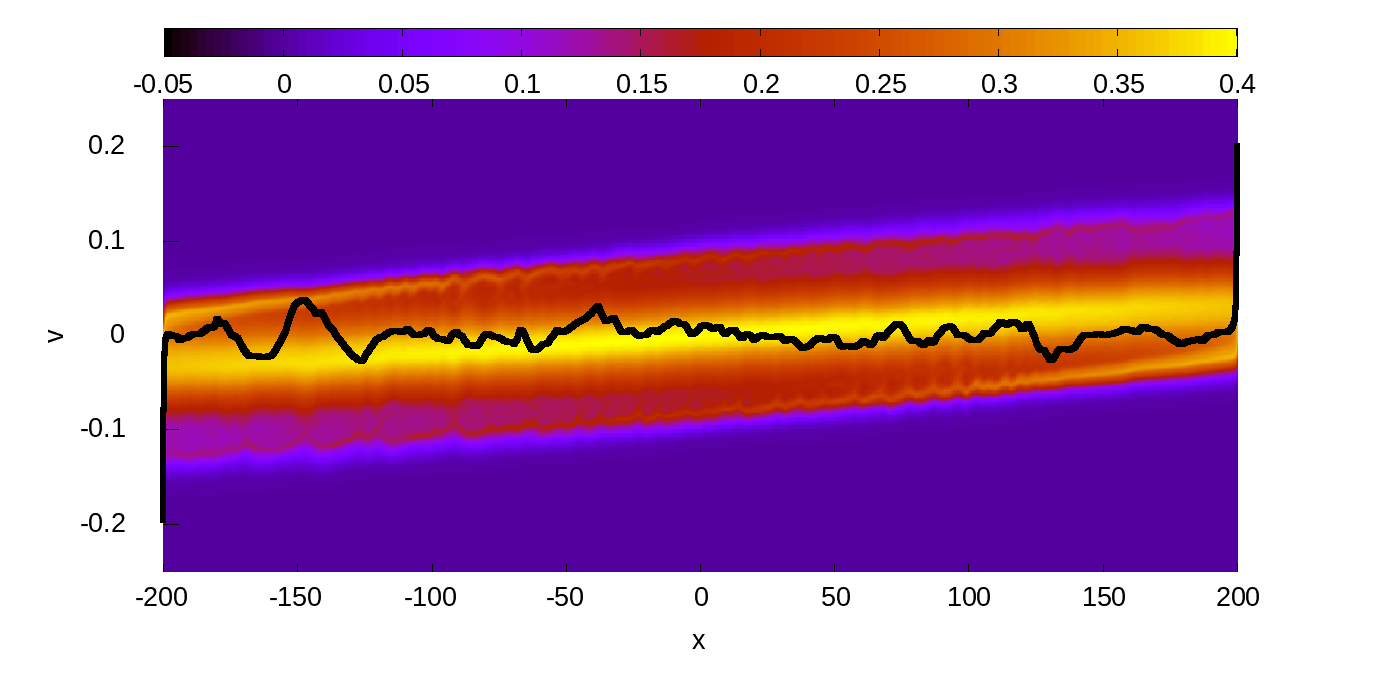}
\caption{The electron particle density for the continuously injected plasma at time $t=7200$ is shown for different limiters. Top left: reference solution, $6000 \times 1200$ degrees of freedom, meanerr indicator and line WENO modifier. Top right: minmod indicator and line WENO modifier, $1200 \times 600$ degrees of freedom. Bottom left: meanerr indicator and line WENO modifier, $1200 \times 600$ degrees of freedom. Bottom right: the newly proposed sLdG limiter,  $1200 \times 600$ degrees of freedom. In all simulations a refinement ratio of $8$ is used close to the wall. \label{fig:density-conting}}
\end{figure}

To investigate the behavior of the limiters in more detail, let us consider Figure \ref{fig:density-conting}, where we show the density function and the electric field at time $t=7200$. We observe that using the minmod indicator suppresses small scale oscillations, but it creates unphysical large wavelength fluctuations. If the grid is refined those are significantly reduced, but at the cost of increased computational cost. The meanerr based limiters and the newly proposed sLdG limiter perform comparable. Let us also remark that in all cases oscillations at the sharp gradient that occur in the plasma sheath are avoided. 

\subsection{Performance considerations \label{sec:performance}}

In this section we compare the computational performance of the different limiters. All results are achieved on a NVIDIA TitanV GPU. As a baseline we consider running the program without any limiter in the configuration used in the last two sections. In this case the advection in the spatial direction requires $8.5$ seconds and in the velocity direction $9.4$ seconds. All timings reported in this section use $1200 \times 600$ degrees of freedom and the setup of section \ref{sec:blob}. The difference between the spatial and velocity advection is mainly since different spatial cells are adjacent in memory, while the memory access pattern is strided in the velocity direction.

The newly proposed limiter has been designed such that it can directly be computed as part of the main sLdG algorithm. In this case the spatial advection/limiting requires 10.2 seconds and the velocity direction/limiting requires 10.8 seconds. Thus, the run time increases by approximately 15-20\%. This is much cheaper than performing the WENO limiters. Since on GPUs we deal with a memory bound problem, see~\cite{Einkemmer2020a,Einkemmer2022}, reading and writing memory is the most expensive part. Thus, when the limiter is applied as a postprocessing step, additional reads of the entire memory and potential writes (if the cell is modified) have to be performed. The measured timings for the different limiters are reported in table~\ref{tab:timings_weno_limiter}. It can be observed that the meanerr indicator is almost twice as expensive as the minmod indicator. On the other hand, there is no significant performance difference between the simple WENO and the line WENO modifier. 
\begin{table}
\centering
 \begin{tabular}{l|S[table-format=2.1] S[table-format=2.1] S[table-format=2.1] S[table-format=2.1]}
 & \text{meanerr+simple} & \text{meanerr+line} & \text{minmod+simple} & \text{minmod+line} \\ \hline
 x-direction: & 10.2s & 10.4s & 4.8s & 5.1s \\
 v-direction: & 6.0s  & 6.0s  & 3.3s & 3.3s \\
 \end{tabular}
\caption{Measured run time for the different WENO limiters. These costs are incurred in addition to the advection step of the sLdG algorithm. }
\label{tab:timings_weno_limiter}
\end{table}

In addition, in the case of mesh refinement close to the wall (as discussed in section \ref{sec:refwall}), the WENO limiters incur additional computational cost due to the required boundary transfer. Even in the two dimensional setting considered here, collecting the data, transferring them from the fine to the coarse grid, and copying them to a buffer such that they can be used is quite expensive. The corresponding run times are collected in table~\ref{tab:bdr_transfer_limiter}. In a case where the different grids are stored on different GPUs or different nodes on a HPC cluster this would be even more severe. For the proposed sLdG limiter this is not an issue as the limiter is directly applied as part of the advection step and thus all the required boundary data is already available.
\begin{table}
\centering
 \begin{tabular}{l|c|c|c}
refinement ratio: & $2$& $4$ & $8$ \\ \hline
time: &  11.2s & 14.8s & 23.1s
 \end{tabular}
\caption{Measured run time for collecting and transferring the boundary data required for the WENO limiter for the spatial mesh refinement close to the wall. These costs are incurred in addition to the advection step of the sLdG algorithm and the actual limiter; the cost of the latter has been reported in table~\ref{tab:timings_weno_limiter}.}
\label{tab:bdr_transfer_limiter}
\end{table}

\section*{Acknowledgements}

This work is supported, in part, by the Austrian Science Fund (FWF) project id: P32143-N32. This project has
received funding, in part, from the European Union’s Horizon 2020 research and innovation programme under the
Marie Skłodowska-Curie grant agreement No 847476. The views and opinions expressed herein do not necessarily
reflect those of the European Commission.

\bibliographystyle{plain}
\bibliography{literature.bib,plasma-physics.bib}

\end{document}